\documentclass[11pt]{article}

\usepackage{amsmath, amssymb, amsthm}
\usepackage{mathrsfs}
\usepackage{xcolor}
\usepackage{comment}

\title{A formula to construct all elements of finite order in Riordan matrix groups}
\author{I. Arévalo, I. Gómez García, A. Luzón, M. M. Martínez-Rico,\\ L. F. Prieto-Martínez}

\date{}

\newtheorem{definition}{Definition}
\newtheorem{remark}[definition]{Remark}

\newtheorem{theorem}[definition]{Theorem}

\newtheorem{coro}[definition]{Corollary}

\newtheorem{lemma}[definition]{Lemma}

\begin{document}

\maketitle

\begin{abstract} We characterize the elements of finite order in the Riordan group over the complex field by providing explicit conditions on their entries. Our main result generalizes a previous characterization of involutions to arbitrary finite orders and gives a constructive procedure to determine all Riordan matrices of a prescribed finite order from a suitable choice of parameters. The proof relies on an inductive analysis of the principal submatrices
%and establishes that every finite-order Riordan matrix is conjugate to a diagonal Riordan matrix. 
We further derive consequences for the Lagrange subgroup and the group of formal power series under composition, obtaining a characterization of finite-order series and a constructive description of the solutions to associated functional equations. 
\end{abstract}

% -------------------------
% INTRODUCTION
% -------------------------

\section{Introduction}

In a group $(G,\cdot)$ with identity element $e$, for some integer $n\geq 2$, we say that a given element $g$ has \textbf{order} $n$ if and only if $n$ is the smallest integer with the property $g^n=e$. The elements of order 2 are called \textbf{involutions}.

We denote by $\mathcal F_0(\mathbb C)$ the set of formal power series $d=d_0+d_1x+d_2x^2+\ldots$ over the complex field, with $d_0\neq 0$. This set is a group with respect to the product of formal power series (\emph{Cauchy product}). In a similar fashion, we denote by $\mathcal F_1(\mathbb C)$ the set of formal power series $h=h_1x+h_2x^2+\ldots$ over the complex field, with $h_1\neq 0$ and, in this case, this set is a group with respect to the composition of formal power series.

The Riordan group is a group of infinite lower triangular matrices and an important tool in combinatorics. It has a rich algebraic structure with its own interest. For more information about this object we refer the reader to the fundational paper \cite{SGWW}, to the surveys \cite{P, S}, and to the books \cite{B, SSBCHMW}. To see some recent examples of the interest in the algebraic study of the Riordan group, and in the elements of finite order in particular, we refer the reader to \cite{CK1, CK2, C, LMP}.

 In the following, we denote by $\mathcal R(\mathbb C)$ the Riordan group over the complex field. For $d\in\mathcal F_0(\mathbb C)$, $h\in \mathcal F_1(\mathbb C)$, we denote by $R(d,h)=[r_{ij}]_{0\leq i,j<\infty}$ the element whose 0-th column contains the coefficients of $d\in \mathcal F_0(\mathbb C)$ and whose $j$-th column contains the coefficients of $dh^j$.
\[
R(d,h)=\begin{array}{c}
\left[
\begin{array}{ccccc}
r_{00} & & & & \\
r_{10} & r_{11} & & & \\
r_{20} & r_{21} & r_{22} & & \\
r_{30} & r_{31} & r_{32} & r_{33} & \\
\vdots & \vdots & \vdots & \vdots & \ddots
\end{array}
\right]
\\[6pt]
\;\begin{array}{ccccc}
\downarrow & \downarrow & \downarrow & \downarrow &  \\
d & dh & dh^2 & dh^3 & \cdots
\end{array}
\end{array}
\]

In this paper, we generalize the main result in \cite{LMP} (Theorem 6 in \cite{LMP} is the case $n=2$ of the following one), providing a formula to construct all elements of finite order in the Riordan group $\mathcal R(\mathbb C)$:

\begin{theorem} \label{theorem.main} For every root of unity of $r_{11}$ of order $n$ and every choice of the elements
$$r_{10},\ldots, r_{n-1,0},r_{n+1,0},\ldots,r_{2n-1,0},\ldots $$
$$r_{21},\ldots, r_{n,1},r_{n+2,1},\ldots,r_{2n,1},\ldots $$

\noindent we can construct a unique Riordan matrix $\mathbf R=[r_{ij}]_{0\leq i,j<\infty}$ of finite order $n$ with 1 in the position $(0,0)$.

Moreover, some $\mathbf R=[r_{ij}]_{0\leq i,j<\infty}\in\mathcal R(\mathbb C)$ such that $r_{00}=1$ and $r_{11}$ is a root of unity of order $n$ is, in turn, an element of order $n$ if and only if

\begin{itemize}

\item For $i\equiv 0\mod n$, $i>0$,
\begin{equation}\label{E:Cond1}
r_{i0}=-\frac{1}{n}\left(\sum_{\substack{i\geq i_1\geq\dots\geq i_{n-1}\geq 0\\ i_j-i_{j+1}<i \; (i=i_0)}}r_{ii_1}\dots r_{i_{n-1}0}\right).    
\end{equation}

\item For $i\equiv 1\mod n$, $i>1$,
\begin{equation}\label{E:Cond2}
 r_{i1}=-\frac{1}{n}\left(\sum_{\substack{i\geq i_1\geq\dots\geq i_{n-1}\geq 1\\ i_j-i_{j+1}<i \; (i=i_0)}}r_{ii_1}\dots r_{i_{n-1}1}\right).   
\end{equation}

\end{itemize}

\end{theorem}

In order to illustrate this result, suppose that we fix $r_{00}=1$, $r_{11}=\mathfrak{i}$  (only for this example, we denote by $\mathfrak{i}$ to be the imaginary unit in $\mathbb C$, to avoid confusion with respect to the index $i$ in the formulas). In this case, $r_{11}$ is a root of unity of order $4$ in $\mathbb C$. We also fix
$$r_{10}=2,\;r_{20}=4,\; r_{30}=8,\; r_{50}=16,\;r_{60}=24,\;r_{70}=32,$$
$$r_{21}=2\mathfrak{i},\;r_{31}=4\mathfrak{i},\; r_{41}=8\mathfrak{i},\; r_{61}=16\mathfrak{i},\;r_{71}=24\mathfrak{i},\;r_{81}=32\mathfrak{i},$$
$$\forall i\geq 9,i\not\equiv 4\mod 4, r_{i0}=0,\forall k\geq 10,k\not\equiv 1\mod 4, r_{k1}=0.$$

\noindent Then we can construct a unique element $[r_{ij}]_{0\leq i,j<\infty}$ of order 4 in $\mathcal R(\mathbb C)$:
\setcounter{MaxMatrixCols}{11}
$$
\begin{bmatrix}
1 & 0 & 0 & 0 & 0 & 0 & 0 & 0 & 0 & 0 & \hdots \\
\boxed{2} & {\mathfrak{i}} & 0 & 0 & 0 & 0 & 0 & 0 & 0 & 0 & \hdots \\
\boxed{4} & \boxed{2\mathfrak{i}} & -1 & 0 & 0 & 0 & 0 & 0 & 0 & 0 & \hdots \\
\boxed{8} & \boxed{4\mathfrak{i}} & -2 & -\mathfrak{i} & 0 & 0 & 0 & 0 & 0 & 0 & \hdots \\
{12} & \boxed{8\mathfrak{i}} & -4 & -2\mathfrak{i} & 1 & 0 & 0 & 0 & 0 & 0 & \hdots \\
\boxed{16} & {12\mathfrak{i}} & -8 & -4\mathfrak{i} & 2 & \mathfrak{i} & 0 & 0 & 0 & 0 & \hdots \\
\boxed{24} & \boxed{16\mathfrak{i}} & -12 & -8\mathfrak{i} & 4 & 2\mathfrak{i} & -1 & 0 & 0 & 0 & \hdots \\
\boxed{32} & \boxed{24\mathfrak{i}} & -16 & -12\mathfrak{i} & 8 & 4\mathfrak{i} & -2 & -\mathfrak{i} & 0 & 0 & \hdots \\
{40} & \boxed{32\mathfrak{i}} & -24 & -16\mathfrak{i} & 12 & 8\mathfrak{i} & -4 & -2\mathfrak{i} & 1 & 0 & \hdots \\
\boxed{0} & {40\mathfrak{i}} & -32 & -24\mathfrak{i} & 16 & 12\mathfrak{i} & -8 & -4\mathfrak{i} & 2 & \mathfrak{i} & \hdots \\
\vdots & \vdots & \vdots & \vdots & \vdots & \vdots & \vdots & \vdots & \vdots & \vdots & \ddots 
\end{bmatrix}.
$$
\setcounter{MaxMatrixCols}{10}

\noindent For this example, the resulting matrix is, in fact, 
$$R(1+2x+4x^2+8x^3+12x^4+16x^5+24x^6+32x^7+40x^8+48x^9,\mathfrak{i}x).$$

\section{Some comments on leading principal submatrices of Riordan matrices}

Given $\mathbf R=[r_{ij}]_{0\leq i,j<\infty}\in\mathcal R(\mathbb C)$ (resp. $\mathbf R=R(d,h)$), let us denote by $\mathbf R_k$ (resp. $R_k(d,h)$) to the leading principal submatrix $[r_{ij}]_{0\leq i,j\leq k}$ \cite{LMMPS, P}. 

We define the set $\mathcal R_k(\mathbb C)=\{\mathbf R_k:\mathbf R\in\mathcal R(\mathbb C)\}$.

Note that, for two different Riordan matrices $\mathbf R,\mathbf S\in\mathcal R(\mathbb C)$, it is possible that $\mathbf R_k=\mathbf S_k$, for some $k$.

Given some $\mathbf R_k\in\mathcal R_k(\mathbb C)$, we denote by $\mathcal E(\mathbf R_k)$ to the subset of $\mathcal R_{k+1}(\mathbb C)$ consisting in those elements obtained from $\mathbf R_k$ adding a row and column. 

\begin{remark} \label{remark.extension} Every element in $ \mathcal E(\mathbf R_k)$ is determined univocally by the choice for the entries in the positions $(k+1,0)$ and $(k+1,1)$. The rest of the entries in $\mathbf R_{k+1}$ are determined by $\mathbf R_k$:
$$\left[\begin{array}{cccc|c}& & & & \\ &  \mathbf R_k & & & \\ & & & & \\ \hline \boxed{r_{k+1,0}}&\boxed{r_{k+1,1}} & r_{k+1,2} & \ldots &r_{k+1,k+1} \end{array}\right].$$

\end{remark}

%%%%%%%%%%%%%%%%%%%%%%%%%%%%%%%%%%%%%%%%%%%%%%%
\section{Proof of Theorem \ref{theorem.main}}

\subsection{Proof of the necessity of the conditions}

Let $\mathbf R=[r_{ij}]_{0\leq i,j<\infty}\in\mathcal R(\mathbb C)$, with $r_{00}=1$, and $r_{11}$ being a root of unity of order $n$.

Let us recall that the entry $(i,j)$ in matrix $\mathbf R^n$ equals to
\begin{equation} \label{eq.formulaproducto}\sum_{i\geq k_1\geq\dots\geq k_{n-1}\geq j}r_{ik_1}\dots r_{k_{n-1}j}. \end{equation}

So, $\mathbf R$ is an element of finite order if and only if its entries satisfy the infinite linear system
\begin{equation}\label{eq.sysinf}
\sum_{i\geq k_1\geq\dots\geq k_{n-1}\geq j}r_{ik_1}\dots r_{k_{n-1}j}=0,  \qquad\text{for }0\leq j<i<\infty. \end{equation}

\begin{itemize}

\item The necessity of Condition \eqref{E:Cond1} is a consequence of System \eqref{eq.sysinf}. If $\mathbf R$ is an element of finite order $n$, then $\mathbf R^n$ is the identity and so, for $i\geq 1$, the $(i,0)$ entry in $\mathbf R^n$ is 0. For  $j=0$ and every $i\geq 1$, note that the coefficient of $m_{i0}$ in the position $(i,0)$ of  $\mathbf R^n$ is 
$$(r_{11}^i)^{(n-1)}+\ldots+r_{11}^i+1=\begin{cases} 0& \text{if } n\not\equiv 0\mod n\\ n&\text{if }n\equiv 0\mod n\end{cases}.$$

\noindent From this last equation one can derive easily Condition \eqref{E:Cond1}.
%$$\begin{cases}=0& \text{if }\\ \neq 0 & \text{in other case.}\end{cases}$$

\item A similar argument shows the neccessity of Condition \eqref{E:Cond2}.

\end{itemize}

\subsection{Proof of the sufficiency of the condition}

Let $\mathbf R\in\mathcal R(\mathbb C)$, with 1 in the position $(0,0)$ and some root of unity of order $n$ in the position $(1,1)$. Then $\mathbf R$ is an element of order $n$ in $\mathcal R(\mathbb C)$ if and only if, for all $k\geq 1$, $\mathbf R_k$ is an element of order $n$ in $\mathcal R_n(\mathbb C)$.

We say that two Riordan matrices $\mathbf R,\mathbf S$ (sim. for $\mathbf R_k,\mathbf S_k$) are \textbf{conjugated}, and we write $\mathbf R\sim\mathbf S$ if there exists a third Riordan matrix $\mathbf P$ such that $\mathbf R=\mathbf P\mathbf S\mathbf P^{-1}$. If $\mathbf R\sim\mathbf D$ (sim. for $\mathbf R_k\sim\mathbf D_k$), where $\mathbf D$ is a diagonal matrix, we say that $\mathbf R$ is \textbf{diagonalizable}.

As happened for the condition for finite order, $\mathbf R$ is conjugated to a diagonal matrix $\mathbf D=R(1,r_{11}x)$ (is diagonalizable), if and only if, for all $k\geq 2$, $\mathbf R_k$ is conjugated to the diagonal matrix $\mathbf D_k=R_k(1,r_{11}x)$ (is diagonalizable).

In the following we prove by induction in $k\geq 1$ the next

\noindent \textbf{Statement (k).} \emph{If $\mathbf R_k$ satisfies Conditions \eqref{E:Cond1} and \eqref{E:Cond2} from Theorem \ref{theorem.main} (obviously, only for $i\leq k$, $i>0$), then   $\mathbf R_{k}$ is an element of order $n$ and  $\mathbf R_k$ is conjugated to $R_k(1,r_{11}x)$. 
}

\noindent If this statement holds for every $k\geq 2$, this implies the sufficiency of the conditions.

For the case $k=1$, note that
$$\mathbf R_1^n=\begin{bmatrix}1\\ r_{10} & r_{11} \end{bmatrix}^n=\mathbf I_1$$

\noindent so $\mathbf R_1$ has order exactly $n$. Moreover, $\mathbf R_1$ is diagonalizable (and so conjugated to $R(1,r_{11}x)$).

Now, for $k\geq 1$, let us assume that $\mathbf R_k$ is an element of finite order $n$ and that it is diagonalizable (so it satisfies the conditions, according to the previous subsection).

 Let us define 
$$\mathcal E_{D}(\mathbf R_k)=\{\mathbf R_{k+1}\in\mathcal E(\mathbf R_k):\mathbf R_{k+1}\sim R_{k+1}(1,r_{11}x) \text{ (is diagonalizable)}\},$$
$$\mathcal E_{FO}(\mathbf R_k)=\{\mathbf R_{k+1}\in\mathcal E(\mathbf R_k):\text{ $\mathbf R_{k+1}$ is of finite order }n\},$$
$$\mathcal E_{C}(\mathbf R_k)=\{\mathbf R_{k+1}\in\mathcal E(\mathbf R_k):\text{ $\mathbf R_{k+1}$ satisfies the conditions}\}.$$

\noindent Observe that $\mathcal E_D(\mathbf R_k)\subset \mathcal E_{FO}(\mathbf R_k)\subset\mathcal E_{C}(\mathbf R_k)$. Moreover, the first set is nonempty (hence the other two are nonempty as well). We only need to prove that $\mathcal E_{C}(\mathbf R_k)\subset \mathcal E_D(\mathbf R_k)$.

%Before further proceeding, let us state the following:

%\begin{itemize}
%    \item \textbf{Conditions to guarantee $\mathcal E_{C}(\mathbf R_{k})\subset \mathcal E_{FO}(\mathbf R_k)$:} To prove  $\mathcal E_{C}(\mathbf R_{k})\subset \mathcal E_{FO}(\mathbf R_k)$, we only need to check that, for any $\mathbf R_{k+1}\in \mathcal E_{C}(\mathbf R_{k})$,
%the last row of $\mathbf R_{k+1}^n$ equals $[0,\ldots, 0,1]$.   Since $\mathbf R_{k+1}^n$ is a submatrix of a Riordan matrix, in fact, the last row must be of the type $[m_{k0},m_{k1},0,\ldots,0,1]$. Moreover, all the entries in the last row of $\mathbf R^n_{k+1}$ except for $m_{k+1,0}$ and $m_{k+1,1}$ are determined by $\mathbf R^n_{k}$. So now we need to study the equations 
%\begin{equation}\label{eq.sys}\begin{cases}\displaystyle{\sum_{k+1\geq k_1\geq\dots\geq k_{n-1}\geq 0}r_{k+1,k_1}\dots r_{k_{n-1}0}=0,}\\\displaystyle{\sum_{k+1\geq k_1\geq\dots\geq k_{n-1}\geq 1}r_{k+1,k_1}\dots r_{k_{n-1}1}=0.} \end{cases}\end{equation}

%\noindent We consider the second of these equations as an equation in the indeterminate $r_{k+1,1}$ and the first one as an equation in the indeterminate $r_{k+1,0},r_{k+1,1}$ (since the rest of the terms in both of them are already determined by $\mathbf R_k$). 

%\item 

Let us remark that, provided that  $\mathbf R_k=\mathbf P_k\mathbf D_k\mathbf P_k^{-1}$, to check that some $\mathbf R_{k+1}=[r_{ij}]_{0\leq i,j\leq k+1}\in\mathcal E(\mathbf R_{k})$ also lies in $\mathcal E_D(\mathbf R_k)$, we need to decide if there exists some $\mathbf P_{k+1}=[p_{ij}]_{0\leq i,j\leq k+1}\in\mathcal E(\mathbf P_k)$ such that $\mathbf R_{k+1}=\mathbf P_{k+1}\mathbf D_{k+1}\mathbf P_{k+1}^{-1}$. This is equivalent, in turn, to check only the equality for the positions $(k+1,0),(k+1,1)$, that is, we need to study the following system of equations in the indeterminates $p_{k+1,0},p_{k+1,1}$
$$\begin{cases}p_{k+1,0}=r_{k+1,0}p_{00}+\ldots+r_{k+1,k+1}p_{k+1,0}\\  \omega_1p_{k+1,1}=r_{k+1,1}p_{11}+\ldots+r_{k+1,k+1}p_{k+1,1}\end{cases}$$

\noindent which is, in turn, equivalent to
\begin{equation}\label{eq.sys2}
\begin{cases}(1-r_{k+1,k+1})p_{k+1,0}=r_{k+1,0}p_{00}+\ldots+r_{k+1,k}p_{k0}\\  (r_{11}-r_{k+1,k+1})p_{k+1,1}=r_{k+1,1}p_{11}+\ldots+r_{k+1,k}p_{k1}.\end{cases}
\end{equation}

%\noindent We consider the first equation as an equation in the intdeterminate $p_{k0}$ and the second equation as an equation in the indeterminate $p_{k1}$. So System \eqref{eq.sys2} is a linear system in the indeterminates $p_{k0},p_{k1}$.

%\end{itemize}

We have three cases:

\begin{itemize}

\item Case $k+1\not\equiv 0,1\mod n$. Note that, in this case, $\mathcal E_C(\mathbf R_k)=\mathcal E(\mathbf R_k)$.

%We first prove that every element in $\mathcal E(\mathbf R_k)$ is an element of finite order.

%We have already discussed that the coefficient of $r_{k+1,1}$ in the second equation is 0. So either there is no element of finite order $n$ or, for every $r_{k+1,1}$, the second equation in the system holds.  But since the set of elements of finite order in $\mathcal E(\mathbf R_k)$ is nonempty, then, for every $r_{k+1,1}$, the second equation in the system holds. 

%Once we have fixed the value of $r_{k+1}$, we repeat the same argument to prove that, for every $r_{k+1,0}$, the first equation holds.

%Now we prove that, 
For every $\mathbf R_{k+1}=[r_{ij}]_{0\leq i,j\leq k}\in\mathcal E(\mathbf R_k)$, the first equation in System \eqref{eq.sys2} has a unique solution in the indeterminate $p_{k+1,0}$ and, similarly, the second equation has a unique solution in the indeterminate $p_{k+1,1}$.

%, every element in $\mathcal E(\mathbf R_k)$ is of the type $\mathbf P_{k+1}\mathbf D_{k+1}\mathbf P_{k+1}^{-1}$, for some $\mathbf P_{k+1}\in\mathcal E(\mathbf P_k)$. F

%from System \eqref{eq.sys2}, considering the first equation as an equation in the indeterminate $p_{k0}$ (the rest of terms are known), we can see that this equation has a unique solution. Similarly, considering the second equation as an equation in the indeterminate $p_{k1}$, we can see that this equation has a unique solution too.

\item Case $k+1\equiv 0\mod n$. There is a unique element in $\mathcal E_C(\mathbf R_k)$ for each choice of the entry $(k+1,1)$.

%Now we prove that, provided that  $\mathbf R_k=\mathbf P_k\mathbf D_k\mathbf P_k^{-1}$, every element in $\mathcal E(\mathbf R_k)$ is of the type $\mathbf P_{k+1}\mathbf D_{k+1}\mathbf P_{k+1}^{-1}$, for some $\mathbf P_{k+1}=[p_{ij}]_{0\leq i,j\leq k+1}\in\mathcal E(\mathbf P_k)$. The position $(k,1)$ in $\mathbf P_{k+1}\mathbf D_{k+1}\mathbf P_{k+1}^{-1}$ is of the type

For every choice of the parameter $r_{k+1,1}$, the second equation in System \eqref{eq.sys2} has a unique solution in the indeterminate $p_{k+1,1}$. But in this case, since the coefficient of $p_{k+1,0}$ vanishes, there is only one choice of the parameter $r_{k+1,0}$ that makes the first equation in System \eqref{eq.sys2} hold.

 Since there is exactly one element in $\mathcal E_D(\mathbf R_k)$ for each choice of $r_{k+1,1}$, there is exactly one element in $\mathcal E_C(\mathbf R_k)$ for each choice of $r_{k+1,1}$ and $\mathcal E_D(\mathbf R_k)\subset \mathcal E_C(\mathbf R_k)$, we conclude that $\mathcal E_D(\mathbf R_k)= \mathcal E_C(\mathbf R_k)$.

\item Case $k+1\equiv 1\mod n$. There is a unique element in $\mathcal E_C(\mathbf R_k)$ for each choice of the entry $(k+1,0)$. This case is very similar to the previous one. We omit details.

\end{itemize}

%%%%%%%%%%%%%%%%%%%%%%%%%%%%%%%%%%%%%%%%%%%%%%%%%%%%%%%%%%%%%%%%%%%%%%%%%%%%%%%%%%%%
\section{Consequences of Theorem \ref{theorem.main}}

\subsection{General elements of finite order}

Let us start with the following:

\begin{lemma} If $\mathbf R=[r_{ij}]_{0\leq i,j<\infty}$ is some element of finite order $n$ such that $r_{00}=1$, then $r_{11}$ needs to be a root of unity or order (exactly, not less than) $n$.
    
\end{lemma}

\begin{proof}
    Suppose that $r_{11}$ is or order $m$, for some proper divisor $m$ of $n$. Then $\mathbf R^{n/m}$ is some element of finite order with ones in the main diagonal. A simple proof by induction shows that the only Riordan matrix $\mathbf S$ such that $\mathbf S^{n/m}=R(1,x)$ is $\mathbf S=R(1,x)$. We omit details.
\end{proof}
    
Moreover, any element of in $\mathcal R(\mathbb C)$ can be obtain as a product $\mathbf D\mathbf R$, for $\mathbf D$ being a diagonal Riordan matrix with constant diagonal and $\mathbf R$ having a 1 in position $(0,0)$.

As a consequence, it is immediate to extend Theorem \ref{theorem.main} to a formula to construct any element $\mathbf R$ of order $n$ in $\mathcal R(\mathbb C)$.

\subsection{General elements of finite order and remarks on conjugacy}

As a consequence of our proof of Theorem \ref{theorem.main} (and of the argument in the previous subsection) we also recover the following known theorem: 

\begin{theorem}[Statement A in Theorem 2 in \cite{C}, Conjugacy Theorem] \label{theorem.cohen} A Riordan matrix  $\mathbf R=[r_{ij}]_{0\leq i,j<\infty}$ is an element of order $n$ if and only if $\mathbf R\sim R(r_{00},r_{11}x)$.

\end{theorem}

This result is closely related to a classical linearization theorem, which states that every formal power series $\omega=\omega_1x+\omega_2x^2+\ldots$ of finite order $n$ in $\mathcal F_1(\mathbb C)$ is conjugate to its linear part $\omega_1x$. We refer the reader to \cite{CohenA} and the references therein for further details concerning this theorem.

Moreover, as a direct consequence of the argument in the proof of Theorem \ref{theorem.main}, we can also prove the following (we omit details):

\begin{coro} Let $\mathbf R$ be an element of order $n$ in $\mathcal R(\mathbb C)$. For every choice of the entries 
$$p_{00},p_{n0},p_{2n,0},\ldots$$
$$p_{11},p_{n+1,1},p_{2n+1,1},\ldots$$

\noindent we can construct a unique Riordan matrix in $\mathcal R(\mathbb C)$ such that $\mathbf R=\mathbf P\mathbf D\mathbf P^{-1}$, where $\mathbf D$ is the diagonal matrix whose elements coincide with those in the diagonal of $\mathbf R$.

\end{coro}

\subsection{Comments on the Lagrange subgroup and $\mathcal F_1(\mathbb C)$}

The \textbf{Lagrange subgroup} is the subgroup $\mathcal L(\mathbb C)<\mathcal R(\mathbb C)$ consisting in those Riordan matrices of the type $R(1,h)$.

\begin{coro} Some element $\mathbf B=[b_{ij}]_{0\leq i,j<\infty}$ in the {Lagrange subgroup} is of finite order $n$ if and only if $b_{11}$ is of finite order $n$ and, for $i\equiv 1\mod n$,
$$b_{i1}=-\frac{1}{n}\left(\sum_{\substack{i\geq i_1\geq\dots\geq i_{n-1}\geq 1\\ i_j-i_{j+1}<i \; (i=i_0)}}b_{ii_1}\dots b_{i_{n-1}1}\right).$$

\end{coro}

\begin{proof} From Condition \eqref{E:Cond1} in Theorem \ref{theorem.main}, we can see that, for every $i$ such that $i\equiv 0\mod n$ and $i\geq 1$,
$$\forall k, \;1\leq k< i, \;b_{k0}=0\Rightarrow b_{i0}=0.$$
\end{proof}

Since $\mathcal L(\mathbb C)$ and $\mathcal F_1(\mathbb C)$ are isomorphic, then the previous corollary also provide a characterization of the elements of order $n$ in $\mathcal F_1(\mathbb C)$, which is a subject of its own interest (we refer the reader to \cite{BB} for more about this topic).

\begin{coro} \label{coro.fenix} Some element $\omega=\omega_1x+\omega_2x^2+\ldots\in\mathcal F_1(\mathbb C)$ is of finite order if and only if:
\begin{itemize}

\item The coefficient $\omega_{1}$ is of order $n$ in the multiplicative group $\mathbb C^*$.

\item  For $i\equiv 1\mod n$,
$$\omega_{i}=-\frac{1}{n}\left(\sum_{\substack{i\geq i_1\geq\dots\geq i_{n-1}\geq 1\\ i_j-i_{j+1}<i \; (i=i_0)}}b_{ii_1}\dots b_{i_{n-1}1}\right)$$

\noindent for $b_{ij}=[x^i]\omega^j$ (or, equivalently, $R(1,\omega)=[b_{ij}]_{0\leq i,j<\infty}$).

\end{itemize}

\end{coro}

\subsection{Consequences related to functional equations}
\begin{remark}\label{remark.functional}
The fact of $R(\alpha,\omega)$ being an element of order $n$ implies $(R(\alpha,\omega))^n=\mathbf I$ and this is equivalent to $\alpha,\omega$ satisfying the functional equations:
$$\alpha\alpha(\omega)\ldots\alpha(\omega^{\circ (n-1)})=1,\qquad \omega^{n-1}=x. $$
\end{remark}

\noindent So we have the following:

\begin{coro} Let $\omega$ be an element of order $n$ in $\mathcal F_1(\mathbb C)$. For every choice of the coefficients 
$$\alpha_2,\ldots\alpha_{k-1},\alpha_{k+1},\ldots, \alpha_{2k-1},\ldots$$
\noindent we can construct a unique solution $\alpha=\alpha_0+\alpha_1x+\alpha_2x^2+\ldots$ for the functional equation 
$$\alpha\alpha(\omega)\ldots\alpha(\omega^{\circ (n-1)})=1.$$

\end{coro}

Moreover, we believe that Theorem \ref{theorem.main} may have many other applications related to functional equations. As a matter of example, an argument very similar to the one in the proof of Theorem \ref{theorem.main} allow us to show the following statement, that is related to a question proposed by R. Slowik in the conference RART 2026 in Gliwice (concerning some work in progress by her and T. He).

\begin{coro} \label{theorem.roksana} Let $\omega\in\mathcal F_1(\mathbb C)$ be an element of order $n$, for some $n$ even greater than 2. The following statements are equivalent:

\noindent (a) $g=g_0+g_1x+\ldots\in\mathcal F_0(\mathbb C)$, with $g_0=1$, is a solution of $gg(\omega)=1$.

\noindent (b) There exists some $\alpha\in\mathcal F_0(\mathbb C)$ such that $R(\alpha,\omega)$ is an element of order $n$ and $g=\alpha\alpha(\omega\circ \omega)\ldots\alpha(\omega^{\circ (n-2)})  $.

\end{coro}

\begin{proof}
    Statement (b) clearly implies Statement (a). To prove that Statement (a) implies Statement (b) we can do an argument, similar to the one in the proof of Theorem \ref{theorem.main}. We omit details, but we include a sketch for the proof:

    \begin{itemize}
\item     For every formal power series $d=d_0+d_1x+\ldots$, let us denote 
$$\mathbf v_d=[d_0,d_1,\ldots]^T,\qquad \mathbf v_{d,k}=[d_0,\ldots, d_k]^T.$$

        \item The \emph{statement} is, for each $k\geq 1$, 
    $$R_k(g,\omega)\mathbf v_{g,k}=\mathbf v_{1,k}\text{ if and only if }\mathbf v_{g,k}=R_k(\alpha,\omega)^{n-1}\mathbf v_{\alpha,k}.$$

\item We denote by $\mathbb C_{\leq k}[x]$ the set of polynomials of degree less than or equal to $k$. For each $g=g_0+\ldots+g_kx^k\in\mathbb C_{\leq k}[x]$, we consider the sets
$$\mathcal E(g)=\{\widetilde g\in\mathbb C_{\leq k+1}[x]:\widetilde g=g_0+\ldots+g_kx^k+\widetilde g_{k+1}x^{k+1}\},$$
$$\mathcal E_C(g)=\{\widetilde g\in\mathcal E(g):R_{k+1}(\widetilde g,\omega)\mathbf v_{\widetilde g,k+1}=\mathbf v_{1,k+1}\},$$
$$\mathcal E_P(g)=\{\widetilde g\in\mathcal E(g):\mathbf v_{\widetilde g,k+1}=R_{k+1}(\alpha,\omega)^{n-1}\mathbf v_{\alpha,k+1},\text{ for some }\alpha\in \mathbb C_{\leq k+1}[x]\}.$$

We know that, for any polynomial $g$, $\mathcal E_P(g)\subset \mathcal E_C(g)$ and both sets are nonempty. 

\item For $g\in\mathbb C_{\leq n}$, we need to study the equation 
$$R_{k+1}(\widetilde g,\omega)\mathbf v_{\widetilde g,k+1}=\mathbf v_{1,k+1}$$

\noindent where everything is determined except for the term of greatest order in $\widetilde g$. 

We can see that, for $k+1\not\equiv n/2\mod n$, $\mathcal E_C(g)$ has exactly 1 element. Since $\mathcal E_C(g)$ is nonempty, in the case $k+1\equiv n/2$, $\mathcal E_C(g)$ has infinitely many elements. 

\item Finally, for $g\in\mathbb C_{\leq n}$, we need to proof that each element in $\mathcal E_C(g)$, is also in $\mathcal E_P(g)$. To do so, we study the existence of solutions of the linear equation corresponding to the last position in 
$$\mathbf v_{\widetilde g,k+1}=R_{k+1}(\alpha,\omega)^{n-1}\mathbf v_{\alpha,k+1}.$$

\noindent We suppose that $\alpha=\alpha_0+\ldots+\alpha_{k+1}x^{k+1}\in\mathbb C_{\leq k+1}[x]$, and that the only indeterminate is $\widetilde \alpha_{k+1}$. In this case by case study we use Theorem \ref{theorem.main}).

    \end{itemize}

\end{proof}

\subsection{Elements of finite order out of other elements of finite order and some comments on the A-sequence}

Let $\mathbf R=R(\alpha,\omega)=[r_{ij}]_{0\leq i,j<\infty}$ be some element of finite order $n$ in $\mathcal R(\mathbb C)$. 

Since $\mathbf R$ is lower triangular, the matrix $\widetilde{\mathbf R}=R\left(\alpha\frac{ \omega}{x},\omega\right)$, that is, the one obtained from $\mathbf R$ by deleting the 0-th column and row, is again an element of order $n$ in $\mathcal R(\mathbf C)$.

Similarly, for every $k\in\mathbb Z$, the matrix $R\left(\alpha\left(\frac{ \omega}{x}\right)^k,\omega\right)$, obtained from $\mathbf R$ by deleting $k$ rows and $k$ columns, is again an element of order $n$ in $\mathbf R(\mathbf C)$.

Also, for every $k_1,k_2\in\mathbb N$, every submatrix of the type $[r_{i,j}]_{k_1\leq i,j\leq k_2}$ is an element of finite order $n$ in $\mathcal R_{k+2-k_1}(\mathbb C)$.

The Riordan group has a well-known construction pattern involving the so-called A-sequence. For a given Riordan matrix \(R(d,h)\), the generating function \(A\) of the A-sequence is defined by
\[
h^{-1}=\frac{x}{A}\qquad \text{or, equivalently,}\qquad A=\frac{x}{h^{-1}}.
\]

\noindent The results presented in this paper have several consequences for the generating function of the A-sequence of elements of finite order.

As a first example, we note that a Riordan matrix \(R(1,h)\) is of order \(n\) if and only if \(R(1,\frac{x}{A})\) is also of order \(n\), where \(A\) denotes the generating function of the A-sequence of \(R(1,h)\). Thus, the generating functions \(A\) associated with elements of finite order \(n\) follow a construction pattern similar to that of the elements \(\omega\) of finite order \(n\) in \(\mathcal F_1(\mathbb C)\) (Corollary \ref{coro.fenix}).

If \(h\) is an element of order \(n=2\) in \(\mathcal F(\mathbb C)\), then \(h^{-1}=h\). Hence, the generating function of the A-sequence satisfies \(h=\frac{x}{A}\). Similar generalizations can be derived for finite orders \(n>2\).

We hope to explore these questions in future research.

% -------------------------
% BIBLIOGRAPHY
% -------------------------

\end{document}